\documentclass[11pt]{article}
\usepackage[margin=1in]{geometry}
\usepackage{amsthm,amsmath}
\usepackage{latexsym,amssymb,amsmath}
\usepackage{graphicx,float,color,fancybox,shapepar,setspace,hyperref}
\usepackage{subfigure}
\usepackage{pgf,tikz}
\usepackage[affil-it]{authblk}
\usepackage{indentfirst}
\usepackage{hyperref}
\usepackage[section]{placeins}
\hypersetup
{
colorlinks=true,
linkcolor=blue,
filecolor=blue,
urlcolor=blue,
citecolor=cyan,
}
\usetikzlibrary{arrows}
\newtheorem{thm}{Theorem}[section]

\newtheorem{lem}[thm]{Lemma}
\newtheorem{prop}[thm]{Proposition}
\newtheorem{defn}[thm]{Definition}

 \def\qed{\hfill\square}

\def\~{\sim}
\def\l{\ell}
\def\qed{ \hfill $\blacksquare$}
\def\pf{\medskip\noindent {\emph{Proof}.}~}

\begin{document}

\title{Connectivity of friends-and-strangers graphs on random pairs}

\author{Lanchao Wang,  Yaojun Chen\thanks{Corresponding author. 
 Email: 
\href{mailto://yaojunc@nju.edu.cn}{ yaojunc@nju.edu.cn}.
}}
 \affil{ { \small {Department of Mathematics, Nanjing University, Nanjing 210093, China}}
 }
\maketitle

\begin{abstract}
Consider two graphs $X$ and $Y$, each with $n$ vertices. The friends-and-strangers graph $\mathsf{FS}(X,Y)$ of $X$ and $Y$ is a graph with vertex set consisting of all bijections $\sigma :V(X) \mapsto V(Y)$, where two bijections $\sigma$, $\sigma'$ are adjacent if and only if they differ precisely on two adjacent vertices of $X$, and the corresponding mappings are adjacent in $Y$. The most fundamental question that one can ask about these friends-and-strangers graphs is whether or not they are connected. Alon, Defant, and Kravitz showed that if $X$ and $Y$  are two independent  random graphs in $\mathcal{G}(n,p)$,  then the threshold probability guaranteeing the connectedness of $\mathsf{FS}(X,Y)$ is $p_0=n^{-1/2+o(1)}$, and suggested to investigate the general asymmetric situation, that is,    
$X\in \mathcal{G}(n,p_1)$ and $Y\in \mathcal{G}(n,p_2)$. In this paper, we show that if  $p_1 p_2 \ge p_0^2=n^{-1+o(1)}$ and $p_1, p_2 \ge  w(n) p_0$, where $w(n)\rightarrow 0$  as $n\rightarrow \infty$,  then $\mathsf{FS}(X,Y)$ is connected with high probability, which extends the result  on $p_1=p_2=p$, due to  Alon, Defant, and Kravitz.


\end{abstract}

\section{Introduction} 
All graphs considered in this paper are finite and simple. The vertex set and edge set of a graph $G$ are denoted by $V(G)$ and $E(G)$, respectively. For any $S\subseteq V(G)$, $G|_S$ denotes the induced subgraph of $G$ by $S$.
For a vertex $v\in V(G)$, the degree $d(v)$ of $v$ is the number of edges incident with $v$ in $G$ and the maximum degree of $G$ is $\Delta(G)=\max\{d(v)~|~v\in V(G)\}$.  Let $K_n, S_n, P_n$ and $C_n $ denote a complete graph, a star, a path and a cycle of order $n$, respectively. A complete bipartite graph on $s+t$ vertices is denoted by $K_{s,t}$. Let $[n]=\{1,2,...,n\}$.
The friends-and-strangers graphs were introduced by Defant and Kravitz \cite{DK}, which are defined as follows. 

 \begin{defn} (Defant and Kravitz \cite{DK})
Let $X$ and $Y$ be two graphs, each with $n$ vertices. The friends-and-strangers graph $\mathsf{FS}(X,Y)$ of $X$ and $Y$ is a graph with vertex set consisting of all bijections from $V(X)$ to $V(Y)$, two such bijections $\sigma$, $\sigma'$ are adjacent if and only if they differ precisely on two adjacent vertices, say $a,b\in X$ with $\{a,b\}\in E(X)$, and the corresponding mappings are adjacent in $Y$, i.e.
\begin{itemize}
\item  $\{\sigma(a), \sigma(b)\} \in E(Y)$;
\item  $\sigma(a)=\sigma'(b)$, $\sigma(b)=\sigma'(a)$ and $\sigma(c)=\sigma'(c)$ for all $c\in V(X)\backslash \{a,b\}.$
\end{itemize}
\end{defn}
When this is the case, we refer to the operation that transforms $\sigma$ into $\sigma'$ as an $(X,Y)$-{\it friendly swap}.

 The friends-and-strangers graph $\mathsf{FS}(X,Y)$ gets its name from the following interpretation.  Corresponds $n$ people to $V(Y)$ and $n$ positions to $V(X)$. Two people are friends if and only if their corresponding vertices are adjacent in $Y$.  Two positions are adjacent if and only if their corresponding vertices are adjacent in $X$. Suppose that these $n$ people stand on these $n$ positions such that each person stands on precisely one position, which corresponds to a bijection from $V(X)$ to $V(Y)$. At any point of time, two people can swap their positions if and only if they are friends and the two positions they stand are adjacent. An immediate question is  how various configurations can be reached from other configurations when  multiple such swaps are allowed. This is precisely the information that is encoded in $\mathsf{FS}(X,Y)$. Note that the components of $\mathsf{FS}(X,Y)$ correspond to the equivalence classes of mutually-reachable (by the multiple swaps described above) configurations, so the connectivity, including the giant component phenomenon and so on, is the basic aspect  of interest in  friends-and-strangers graphs. 

 A well-known example is the $15$-puzzle. The numbers from $1$ to $15$ are placed on a $4\times 4$ grid. At each time, two numbers are forbidden to swap their positions, while the ``empty" is allowed to swap its position with any number whose positions are adjacent to its position. By the interpretation above, this game corresponds to $\mathsf{FS}(G_{4\times4},S_{16})$, where $G_{4\times4}$ is the $4\times4$ grid graph and the center of the $S_{16}$ associates with the ``empty'' in the game.
Wilson \cite{W} generalized the $15$-puzzle to $\mathsf{FS}(X,S_{n})$. Friends-and-strangers graphs also generalize  many other objects, for example, $\mathsf{FS}(K_n,Y)$ is isomorphic to the Cayley graph of $\mathfrak{S_n}$ generated by the transpositions corresponding to $E(Y)$, where $\mathfrak{S_n}$ denotes the symmetric group consisting of all permutations of the numbers $1,\dots,n$. 

Let $\mathcal{G}(n,p)$ denote the  Erd\H{o}s-R\'{e}nyi random graphs with $n$ vertices and edge-chosen probability $p$. A sequence of events $\{ \mathcal{A}_n \}$ is said to occur {\it with high probability}, abbreviated to {\it w.h.p.}, if $\Pr (\mathcal{A}_n) \to 1$ as $n\to \infty$.  Moreover, for two functions $f(n)$ and $g(n)$ of $n$, we use  ``$f(n)\gg g(n)$'' instead of  ``$g(n)=o(f(n))$'' throughout this paper. 

The questions and results in literature on the friends-and-strangers graph $\mathsf{FS}(X,Y)$ roughly fall in three types: One or both of $X,Y$ are concrete graphs,  $X$ and $Y$ have an extremal structure, or both $X$ and $Y$ are random ones from  $\mathcal{G}(n,p)$.
In the case when  both $X$ and $Y$ are not random graphs,
Defant and Kravitz \cite{DK}  derived many  basic properties of $\mathsf{FS}(X,Y)$, the structure of $\mathsf{FS}(P_n,Y)$ and $\mathsf{FS}(C_n,Y)$, and general conditions for $\mathsf{FS}(X,Y)$ to be connected. Alon, Defant, and Kravitz \cite{ADK} studied the minimum degree condition for both $X$ and $Y$ that guarantees the connectedness of $\mathsf{FS}(X,Y)$. Bangachev \cite{B} studied the generalized versions of two problems in \cite{ADK} concerning the conditions on minimum degree. Jeong  studied the diameter (the largest distance between any two vertices) and the girth (the length of a shortest cycle) of a friends-and-strangers graph in \cite{J1}  and \cite{J2}, respectively. Defant, Dong, Lee and Wei  \cite{DDLW} studied some new general conditions for $\mathsf{FS}(X,Y)$ to be connected and the cycle spaces of $\mathsf{FS}(C_n,Y)$. Milojevi\'c introduced a generalization of friends-and-strangers graphs in which vertices of the starting graphs are allowed to have multiplicities.
When both $X$ and $Y$ are random graphs from $\mathcal{G}(n,p)$, the most interesting problem maybe is  the threshold for the probability $p$ at which $\mathsf{FS}(X,Y)$ changes from disconnected w.h.p. to connected w.h.p.
Alon, Defant, and Kravitz \cite{ADK} studied this problem for when both $X$ and $Y$ are random graphs in $\mathcal{G}(n,p)$ or random bipartite graphs in $\mathcal{G}(K_{r,r},p)$. They showed that if $X$ and $Y$ are random graphs in $\mathcal{G}(n,p)$, then the threshold probability guaranteeing the connectedness of $\mathsf{FS}(X,Y)$ is $p=n^{-1/2+o(1)}$. More precisely, they proved the following.

%

\begin{thm} (Alon, Defant, and Kravitz \cite{ADK}) \label{q}
Fix some small $\varepsilon >0$, and let $X$ and $Y$ be two graphs independently chosen from $\mathcal{G}(n,p)$, where $p=p(n)$ depends on $n$. If 
\[
\begin{split}
p\le \frac{2^{-1/2}-\varepsilon}{n^{1/2}},
\end{split}
\]
then $\mathsf{FS}(X,Y)$ is disconnected w.h.p. If
\[
\begin{split}
p\ge \frac{\exp(2(\log n)^{2/3})}{n^{1/2}},
\end{split}
\]
then $\mathsf{FS}(X,Y)$ is connected w.h.p.
\end{thm} 

It is clear that both $X$ and $Y$ are taken from $\mathcal{G}(n,p)$ in Theorem \ref{q}, that is, $X$ and $Y$ are chosen with the same edge-chosen probability $p$. 
However, a more general situation is $X\in\mathcal{G}(n,p_1)$ and $Y\in\mathcal{G}(n,p_2)$. Alon, Defant, and Kravitz \cite{ADK} called the general situation  as asymmetric, and suggested to investigate the connectivity of $\mathsf{FS}(X,Y)$ in this case. In this paper, we focus on discussing the same problem as that in Theorem \ref{q} for the  asymmetric case, and the main result is as below.


\begin{thm} \label{mr}
Fix some small $\varepsilon >0$, and let $X$ and $Y$ be independently chosen random graphs in $\mathcal{G}(n,p_1)$ and $\mathcal{G}(n,p_2)$, respectively, where $p_1=p_1(n)$ and $p_2=p_2(n)$ depend on $n$. Let $p_0=\frac{\exp(2(\log n)^{2/3})}{n^{1/2}}$.  If either

\begin{align*}
p_1 p_2 \le \frac{(1-\varepsilon)/2}{n} \text{ \ and \ }   p_1,p_2 \gg \frac{\log n}{n},
\end{align*}
 \text{or }\begin{align*}
\min\{p_1,p_2\} \le \frac{\log n+c(n)}{n} \text{ for some } c(n) \to -\infty,
\end{align*}
then $\mathsf{FS}(X,Y)$ is disconnected w.h.p.
If

\[
\begin{split}
p_1 p_2 \ge p_0^2 \text{ \ and \ }   p_1, p_2 \ge  \frac{2}{(\log n)^{1/3}} p_0,
\end{split}
\]
then $\mathsf{FS}(X,Y)$ is connected w.h.p.
\end{thm}
\noindent Taking $p_1=p_2=p$ in Theorem \ref{mr},  we get Theorem \ref{q}, and so Theorem \ref{mr} extends Theorem \ref{q} generally.
\section{Preliminaries}

In this section, we list some known results for proving Theorem \ref{mr}. The first two are the basic properties of friends-and-strangers graphs.
\begin{lem} (Defant and Kravitz \cite{DK}) \label{dk} The friends-and-strangers graph
$\mathsf{FS}(X,Y)$ is isomorphic to $ \mathsf{FS}(Y,X)$.\end{lem}
\begin{lem} (Defant and Kravitz \cite{DK}) \label{asd}Let $X,\widetilde{X},Y,\widetilde{Y}$ be graphs on $n$ vertices. If $X$ is a subgraph of $\widetilde{X}$ and $Y$ is a subgraph of $\widetilde{Y}$, then $\mathsf{FS}(X,Y)$ is a subgraph of $\mathsf{FS}(\widetilde{X},\widetilde{Y})$.
\end{lem}

 Note that the friends-and-strangers graph $\mathsf{FS}(K_n,Y)$ is the Cayley graph of the symmetric group $\mathfrak{S_n}$ generated by all transpositions corresponding to $E(Y)$.
\begin{lem} (Godsil and Royle \cite{GR}) \label{ar}
Let $X$ be a graph on $n$ vertices. Then $\mathsf{FS}(K_n,X)$ is connected if and only if $X$ is connected.
\end{lem}

For two graphs $X$ and $Y$  on $n$ vertices,  if $\sigma: V(X)\mapsto V(Y)$ is a graph embedding from $X$ to $\overline{Y}$, the complement of $Y$, then $\sigma$ is an isolated vertex in $\mathsf{FS}(X,Y)$. The following lemma states when such an embedding $\sigma$ exists, that is, when $\mathsf{FS}(X,Y)$ has isolated vertices, which implies $\mathsf{FS}(X,Y)$ is disconnected.

\begin{lem} (Catlin \cite{C}, Sauer and Spencer \cite{SS}) \label{l1} If $X$ and $Y$ are graphs on $n$ vertices satisfying $2\Delta(X)\Delta(Y) < n$, then there exists a bijection $\sigma:V(X)\mapsto V(Y)$ such that for every edge $\{a,b\} $ of $X$, the pair $\{ \sigma(a), \sigma(b) \}$ is not an edge in $Y$.
\end{lem}

For any two vertices $\sigma=\sigma_0$ and $\sigma'=\sigma_\ell$ of $\mathsf{FS}(X,Y)$, if $\sigma$ and $\sigma'$ lie in the same component, 
then there is a path $\sigma_0\sigma_1\cdots\sigma_\ell$ connecting $\sigma$ and $\sigma'$. By the definition of $\mathsf{FS}(X,Y)$, an edge $\{\sigma_i,\sigma_{i+1}\}$ means $\sigma_i$ can be transformed into $\sigma_{i+1}$ by an $(X,Y)$-friendly swap, and so  $\sigma$ can be transformed into $\sigma'$ through a  sequence of
$(X,Y)$-friendly swaps. To consider the connectivity of $\mathsf{FS}(X,Y)$ for $X,Y\in\mathcal{G}(n,p)$,  Alon, Defant, and Kravitz \cite{ADK} introduced the notion of an exchangeable pair of vertices:
 Let $X$ and $Y$ be two graphs on $n$ vertices,   $\sigma:V(X)\mapsto V(Y)$  a bijection and
 $u,v \in V(Y)$. 
 We say $u$ and $v$ are $(X,Y)${\it-exchangeable from $\sigma$} if $\sigma$ and $\tau_{uv} \circ \sigma$ are in the same component, where $\tau_{uv}: V(Y)\mapsto V(Y) $ is the bijection such that $\tau_{uv} (u)=v$, $\tau_{uv} (v)=u$ and $\tau_{uv} (w)=w$ for any $w\in V(Y)\backslash \{u,v\}$. In other words,  we say $u$ and $v$ are $(X,Y)$-exchangeable from $\sigma$ if there  is  a  sequence of
$(X,Y)$-friendly swaps that we can apply to $\sigma$ in order to exchange $u$ and $v$, i.e., there is a  path in $\mathsf{FS}(X,Y)$ that connect $\sigma$ to $\tau_{uv} \circ \sigma$. 
The following lemma gives a sufficient condition for $\mathsf{FS}(X,Y)$ being connected in terms of exchangeable pairs of vertices. 


\begin{lem} (Alon, Defant, and Kravitz \cite{ADK}) \label{ex}
 Let $X, Y$ be two graphs on $n$ vertices, and  $X$ is connected. Suppose for any two vertices $u,v\in Y$ and every $\sigma$ satisfying $\{\sigma^{-1}(u),\sigma^{-1}(v)\} \in E(X)$, the vertices $u$ and $v$ are $(X,Y)$-exchangeable from $\sigma$. Then $\mathsf{FS}(X,Y)$ is connected.
\end{lem}

    In general, it is not easy to know if two vertices $u,v\in V(Y)$ are $(X,Y)$-exchangeable from some bijection $\sigma : V(X)\mapsto V(Y)$. The following lemma provides us a method on how to find a pair of vertices in $Y$, which are  $(X,Y)$-exchangeable from some bijections $\sigma : V(X)\mapsto V(Y)$.

\begin{lem} (Alon, Defant, and Kravitz \cite{ADK}) \label{ex2}
  Let $X$, $Y$ be two graphs on $n$ vertices, and $G$, $H$ be two graphs with vertex set $[m+2]$ such that the vertices $m+1$ and $m+2$ are $(G,H)$-exchangeable from the identity bijection  Id $: [m+2] \mapsto [m+2]$. If there are two graph embeddings $\varphi : [m+2] \mapsto V(X) $ from $G$ to $X$  and $ \psi  : [m+2] \mapsto V(Y)$ from $H$ to $Y$, then the vertices $\psi(m+1)$ and $ \psi(m+2)$ are $(X,Y)$-exchangeable from any bijection $\sigma : V(X)\mapsto V(Y)$ satisfying $\sigma \circ \varphi =  \psi   \circ \text{Id}$.
\end{lem}

\section{Proof of Theorem \ref{mr}}
 
Our main idea for proving Theorem \ref{mr}  comes from \cite{ADK}.
 We will divide the proof of Theorem \ref{mr} into two parts: the disconnected part (Proposition \ref{disconnected}) and the connected part (Proposition \ref{connected}). 
\subsection{Disconnected with high probability}

\begin{prop} \label{disconnected}
Fix some small $\varepsilon >0$, and let $X$ and $Y$ be independently chosen random graphs in $\mathcal{G}(n,p_1)$ and $\mathcal{G}(n,p_2)$, respectively, where $p_1=p_1(n)$ and $p_2=p_2(n)$ depend on $n$. If either 
\begin{align*}
p_1 p_2 \le \frac{(1-\varepsilon)/2}{n} \text{ \ and \ }   p_1,p_2 \gg \frac{\log n}{n},
\end{align*}
 \text{ or }\begin{align*}
\min\{p_1,p_2\} \le \frac{\log n+c(n)}{n} \text{ ~for some } c(n) \to -\infty,
\end{align*}
\\
then  $\mathsf{FS}(X,Y)$ is disconnected w.h.p.
\end{prop}
{\bf\pf}Suppose $p_1 p_2 \le (1-\varepsilon)/(2n)$ and $p_1,p_2 \gg \log n /n$.  Because both $p_1$ and $p_2$ are much larger than $\log{n} /n$, it is well known that the degrees of all vertices in $X$ and $Y$ are $p_1 n(1+o(1))$ and $p_2 n(1+o(1))$ w.h.p., respectively. Consequently, w.h.p., $2\Delta(X)\Delta(Y)=2p_1 p_2 n^2(1+o(1)) \le n(1-\varepsilon)(1+o(1))<n$. By Lemma \ref{l1},  there exists a bijection $\sigma:V(G)\to V(H)$ such that for every edge $\{a,b\} $ of $G$, the pair $\{ \sigma(a), \sigma(b) \}$ is not an edge in $H$ w.h.p., which implies that $\sigma$ is an isolated vertex in $\mathsf{FS}(X,Y)$, and so $\mathsf{FS}(X,Y)$ is disconnected.

 If $\min\{p_1,p_2\} \le (\log n+c(n))/n$, then we may assume $p_1 \le (\log n+c(n))/n$ by Lemma \ref{dk}. Thus, it is well known  $X$ is disconnected w.h.p. since $c(n)\to - \infty$, which is equivalent to the fact that $\mathsf{FS}(X,K_n)$ is disconnected by Lemma \ref{ar}. This  implies that $\mathsf{FS}(X,Y)$ is disconnected w.h.p. since it is a spanning subgraph of $\mathsf{FS}(X,K_n)$ by Lemma \ref{asd}. 
\qed

\subsection{Connected with high probability}
In this section, our main task is to show the following.

\begin{prop} \label{connected}
 Let $X$ and $Y$ be independently chosen random graphs in $\mathcal{G}(n,p_1)$ and $\mathcal{G}(n,p_2)$, respectively. Let  $p_0 =\frac{\exp(2(\log n)^{2/3})}{n^{1/2}}$ and $\l=\frac{(\log n) ^{1/3}}{2}$. If
\[
\begin{split}
p_1 p_2 \ge p_0^2 \text{ ~and~ }  p_1, p_2\ge \frac{1}{{\l}}p_0 ,  
\end{split}
\]
then $\mathsf{FS}(X,Y)$ is connected w.h.p.
\end{prop}

Before starting to prove Proposition \ref{connected}, we need three technical lemmas (\ref{a}-\ref{3}) and some additional notations.


 Let $m$ be a positive integer, $G$ and $H$ be two graphs on vertex set $[m]$, and $\sigma: V(X)\mapsto V(Y)$ be a bijection. Let $V_1,\dots,V_m$ be a list of $m$ pairwise disjoint sets of vertices of $Y$. We say that the pair of graphs $(G,H)$ is {\it embeddable in $(X,Y)$ with respect to the sets $V_1,\dots,V_m$ and the bijection $\sigma$} if there exist vertices $v_i\in V_i$ for all $i\in[m]$ such that for all $i,j\in[m]$, we have
\[
\begin{split}
\{i,j\} \in E(H) \Rightarrow \{v_i,v_j\} \in E(Y)  \text{ and } \\
\{i,j\} \in E(G) \Rightarrow \{ \sigma^{-1}(v_i),\sigma^{-1}(v_j)\} \in E(X).
\end{split}
\]  

 Suppose $q_1,\dots,q_m$ are nonnegative integers satisfying $q_1+\dots+q_m\le n$. We say the pair $(G,H)$ is {\it $(q_1,\dots,q_m)$-embeddable in $(X,Y)$} if the pair $(G,H)$ is embeddable in $(X,Y)$ with respect to every list  $V_1,\dots,V_m$ of pairwise disjoint subsets of $V(Y)$  satisfying $|V_i|=q_i$ for all $i\in[m]$ and every bijection $\sigma: V(X)\mapsto V(Y)$. 
 
 The following technical lemma deals with when we can embed  a pair of small graphs into a pair of large random graphs.

\begin{lem} \label{a}
 Let $m,n,q_1,\dots,q_m$ be positive integers such that $Q=q_1+\dots+q_m\le n$, and  $G$, $H$ be two graphs on the vertex set $[m]$. Let $X$ and $Y$ be independently chosen random graphs in $\mathcal{G}(n,p_1)$ and $\mathcal{G}(n,p_2)$, respectively, where $p_1=p_1(n)$ and $p_2=p_2(n)$ depend on $n$. If for every set $J\subseteq [m]$ satisfying $|E(G|_J)|+|E(H|_J)|\ge 1$ we have 
\[
\begin{split}
{p_1}^{|E(G|_J)|} p_2^{|E(H|_J)|} \prod_{j\in J} q_j \ge 3\cdot 2^{m+1} Q \log{n},
\end{split}
\]
then the probability that the pair $(G,H)$ is $(q_1,\dots,q_m)$-embeddable in $(X,Y)$ is at least $1-n^{-Q}$.
\end{lem}
{\bf\pf}We may assume $|E(G)|+|E(H)|\ge 1$ for otherwise the result is trivial. Fix a list $V_1,\dots,V_m$ of pairwise disjoint subsets of $V(Y)$ satisfying $|V_i|=q_i$ for all $i\in [m]$ and an injection $\iota : \bigcup_{i\in[m]} V_i \mapsto V(X)$. There are at most $n^{q_1}n^{q_2}\cdots n^{q_m} n^Q=n^{2Q}$ ways to make these choices. And then extend $\iota^{-1}$ arbitrarily to a bijection $\sigma : V(X)\mapsto V(Y)$. Note that whether or not $(G,H)$ is embeddable in $(X,Y)$ with respect to the sets $V_1,\dots,V_m$ and the bijection $\sigma$, does not depend on the way in which we extend $\iota^{-1}$ to $\sigma$. We will show that the probability that $(G,H)$ is not embeddable in $(X,Y)$ with respect to the sets $V_1,\dots,V_m$ and the bijection $\sigma$ is at most $n^{-3Q}$. This will imply the desired result by using the union bound of probability.

 Given a tuple $t=(v_1,\dots,v_m)\in V_1\times \dots \times V_m$, let $B_t$ be the ``good'' event that for all $i,j\in[m]$, we have 
\[
\begin{split}
\{i,j\} \in E(H) \Rightarrow \{v_i,v_j\} \in E(Y)  \text{ and } \\
\{i,j\} \in E(G) \Rightarrow \{ \sigma^{-1}(v_i),\sigma^{-1}(v_j)\} \in E(X),
\end{split}
\]  
 i.e. $t=(v_1,\dots,v_m)$ guarantees that $(G,H)$ is embeddable in $(X,Y)$ with respect to the sets $V_1,\dots,V_m$ and the bijection $\sigma$. For tuples $t=(v_1,\dots,v_m)$ and $t'=(v_1',\dots,v_m')$ in $V_1\times \dots \times V_m$ and $J\subseteq [m]$, we write $t\~_J t'$ if $J=\{j\in[m]:v_j=v_j'\}$. We write $t\~t'$ if and only if $t\~_J t'$ for some set $J\subseteq [m]$ satisfying $|E(G|_J)|+|E(H|_J)|\ge 1$. Observe that  if $t \nsim t'$, then the events $B_t$ and $B_{t'}$ are independent. Define
\[
\begin{split}
\Delta=\sum_{t\~t'} \Pr[B_t \land B_{t'}],
\end{split}
\]  
where $B_t \land B_{t'}$ is the event that $B_t$ and $B_{t'}$ both occur and the sum is over all ordered pairs $(t,t')$ such that $t,t' \in V_1\times \dots \times V_m$ and $t\~t'$. 

 Let $\mu $ denote the expected number of the events $B_t$ that occur. We have
\[
\begin{split}
\mu&=\sum_{t\in V_1\times \dots \times V_m} \Pr[B_t]\\
&={p_1}^{|E(G)|} p_2^{|E(H)|} \prod_{j\in [m]} q_j \ge 3\cdot 2^{m+1} Q \log{n},
\end{split}
\]  
where the inequality comes from our hypothesis that $J=[m]$.
  
If $\Delta \le \mu$, then by the Janson Inequality \cite{A}, we have
\begin{align*}
\Pr \left[ \bigwedge\limits_{t\in V_1 \times \dots \times V_m} \overline{B_t}\right] \le e^{-\mu + \Delta/2} \le e^{-\mu/2} \le e^{-(6Q\log n)/2}=n^{-3Q}.
\end{align*}

  If $\Delta \ge \mu$, then by the extended Janson Inequality \cite{A}, and the following claim which gives a lower bound on $\mu^2/ (2\Delta)$, we have
\[
\begin{split}
\Pr \left[ \bigwedge\limits_{t\in V_1 \times \dots \times V_m} \overline{B_t}\right] \le e^{-\mu^2/(2\Delta)} \le e^{-3Q\log n} = n^{-3Q}, \end{split}
\]  as desired.
\vskip 2mm
\noindent {\textbf{Claim.}} $\mu^2/(2\Delta) \ge 3Q\log n$.
\vskip 1mm
{\pf}By the definition of ``$t\~t'$'', we can write
\[
\begin{split}
\Delta=\mathop{\sum_{J\subseteq [m]}}_{|E(G|_J)|+|E(H|_J)|\ge 1} \Delta_J, \text{ \ \ where } \Delta_J=\sum \limits_{t\~_J t'} \Pr[B_t \land B_{t'}].
\end{split}
\]  
 For each $J\subseteq [m]$ with $|E(G|_J)|+|E(H|_J)|\ge 1$, we have
\[
\begin{split}
\Delta_J &\le \left( \prod_{j\in J} q_j\right) p_1^{|E(G|_J)|}p_2^{|E(H|_J)|} \left(\prod_{i\in [m]\backslash J} {q_i}^2\right)p_1^{2(|E(G)|-|E(G|_J)|)}p_2^{2(|E(H)|-|E(H|_J)|)}\\
&=\frac{\mu^2}{p_1^{|E(G|_J)|} p_2^{|E(H|_J)|}\prod_{j\in J}q_j},
\end{split}
\]  
where the inequality comes from the following. The factor $\prod_{j\in J} q_j$ is the number of ways to choose vertices $v_j = v_j'\in V_j$ for all $j\in J$,  and 
$$p_1^{|E(G|_J)|}p_2^{|E(H|_J)|}$$
is the probability that we have

\[
\begin{split}
\{i,j\} \in E(H) \Rightarrow \{v_i,v_j\} \in E(Y)  \text{ and } \\
\{i,j\} \in E(G) \Rightarrow \{ \sigma^{-1}(v_i),\sigma^{-1}(v_j)\} \in E(X)
\end{split}
\]  
for all $i,j\in J$. The factor $\prod_{i\in [m]\backslash J} q_i^2$ is an upper bound on the number of ways to choose the distinct vertices $v_i,v_i'\in V_i$ for all $i\in [m]\backslash J$, and  
$$p_1^{2(|E(G)|-|E(G|_J)|)}p_2^{2(|E(H)|-|E(H|_J)|)}$$
 is the probability that we have 
 \[
\begin{split}
\{i,j\} \in E(H) &\Rightarrow \{v_i,v_j\}, \{v_i',v_j'\} \in E(Y)  \text{ and } \\
\{i,j\} \in E(G) &\Rightarrow \{ \sigma^{-1}(v_i),\sigma^{-1}(v_j)\}, \{ \sigma^{-1}(v_i'),\sigma^{-1}(v_j')\} \in E(X)
\end{split}
\]  
for all $(i,j)\in ([m]\times [m])\backslash (J\times J)$.

 There must be a subset $J^*\subseteq [m]$ such that $|E(G|_{J^*})|+|E(H|_{J^*})|\ge 1$ and $\Delta_{J^*}\ge \Delta/2^m$. Thus we have
\[
\begin{split}
\frac{\mu^2}{2\Delta}\ge \frac{\mu^2}{2^{m+1}\Delta_{J^*}} \ge \frac{\mu^2}{2^{m+1}} \cdot \frac{p_1^{|E(G|_{J^*})|} p_2^{|E(H|_{J^*})|}\prod_{j\in J^*}q_j}{\mu^2}\ge 3Q\log n,
\end{split}
\]  
 where the last inequality comes from our hypothesis that $J=J^*$. \qed \\

 The second technical lemma concerns two sparse graphs $G^*$ and $H^*$ on vertex set $[m+2]$, which was constructed in \cite{ADK}. The definitions of $G^*$ and $H^*$ are as follows.

 Let $n$ be a large integer, $m=\lfloor{(\log n)^{2/3}}\rfloor$ and $\l=\lfloor m^{1/2}/2 \rfloor$. Denote the elements of $[m]$ (written in an arbitrary order) by
\begin{align*}
w,x_1,\dots,x_{\l},y_1,\dots,y_{\l},z_1,\dots,z_{m-2\l-1}.
\end{align*}
 Let $H^{**}$ denote a star $S_m$ on vertex set $[m]$ with center $w$, and  $H^*$ a graph obtained from $H^{**}$ by adding the vertices $m+1$ and $m+2$, along with the additional edges of the form $\{m+1,x_i\}$ and $\{m+2,y_j\}$.

 To describe $G^*$, we first describe a graph $G^{**}$ on vertex set $[m]$, consisting of a cycle $C_m$ with $4$ chords. The vertices of $G^{**}$  are arranged along the $C_m$  in such a way that the vertices $z_1,\dots,z_{12}$ appear in this order when we traverse the $C_m$ anti-clockwise, and the 4 chords are the edges $\{z_1,z_6\}$, $\{z_2,z_4\}$, $\{z_7,z_{12}\}$, $\{z_8,z_{10}\}$. The other vertices of $G^{**}$ lie on the $C_m$ in such a way that the following conditions are satisfied:
\begin{itemize}
\item  The cycle contains the edges $\{z_4,z_5\}$, $\{z_5,z_6\}$, $\{z_{10},z_{11}\}$, $\{z_{11},z_{12}\}$ .
\item  The anti-clockwise distance along the cycle between $z_3$ and $z_5$ is $\l-1$, as is the anti-clockwise distance along the cycle between $z_9$ and $z_{11}$.
\item  The anti-clockwise distance along the cycle from $z_2$ to $z_4$ is even.
\item The $2\l+1$ vertices $w,x_1,\dots,x_{\l},y_1,\dots,y_{\l}$ are placed on the cycle so that the distance in $G^{**}$ between any two of them, as well as the distance in $G^{**}$ between any one of them and any one of the vertices $z_3, z_5, z_9, z_{11}$, is at least $m/(3\l)$.
\item The girth of the graph $G^{**}$ is at least $m/6$.
\end{itemize}

 Because $\l=\lfloor m^{1/2}/2 \rfloor$ and there are exactly $5$ potential small cycles (since there are $4$ chords) in $G^{**}$, the graph $G^{**}$ exists when $n$ is large enough. The graph $G^*$ is obtained from $G^{**}$ by adding the vertices $m+1$ and $m+2$ and the additional edges $\{m+1, m+2\}$, $\{m+1,z_3\}$, $\{m+1,z_{11}\}$, $\{m+2, z_5\}$, $\{m+2,z_9\}$ to $G^{**}$. The graphs $G^{**}$ and $G^*$ are shown in Figure 1,
 where each blue line represents a path of certain length and each of black and red lines represents an edge.
 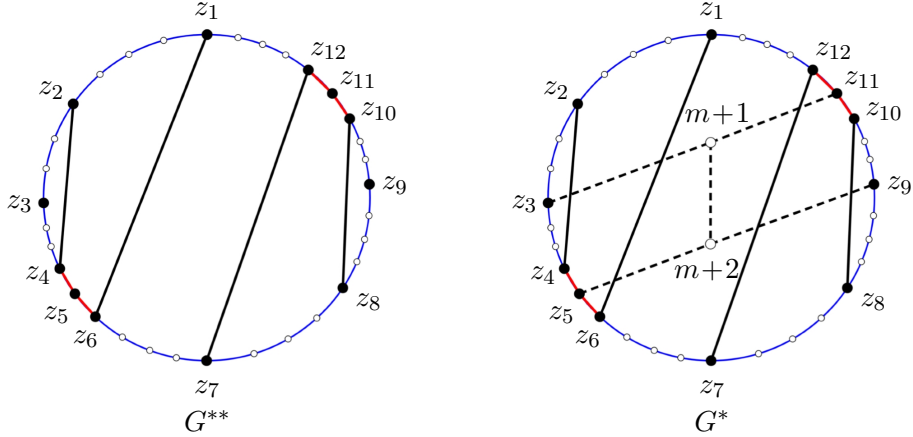
\begin{figure}[h]
\centering
\includegraphics[scale=0.17]{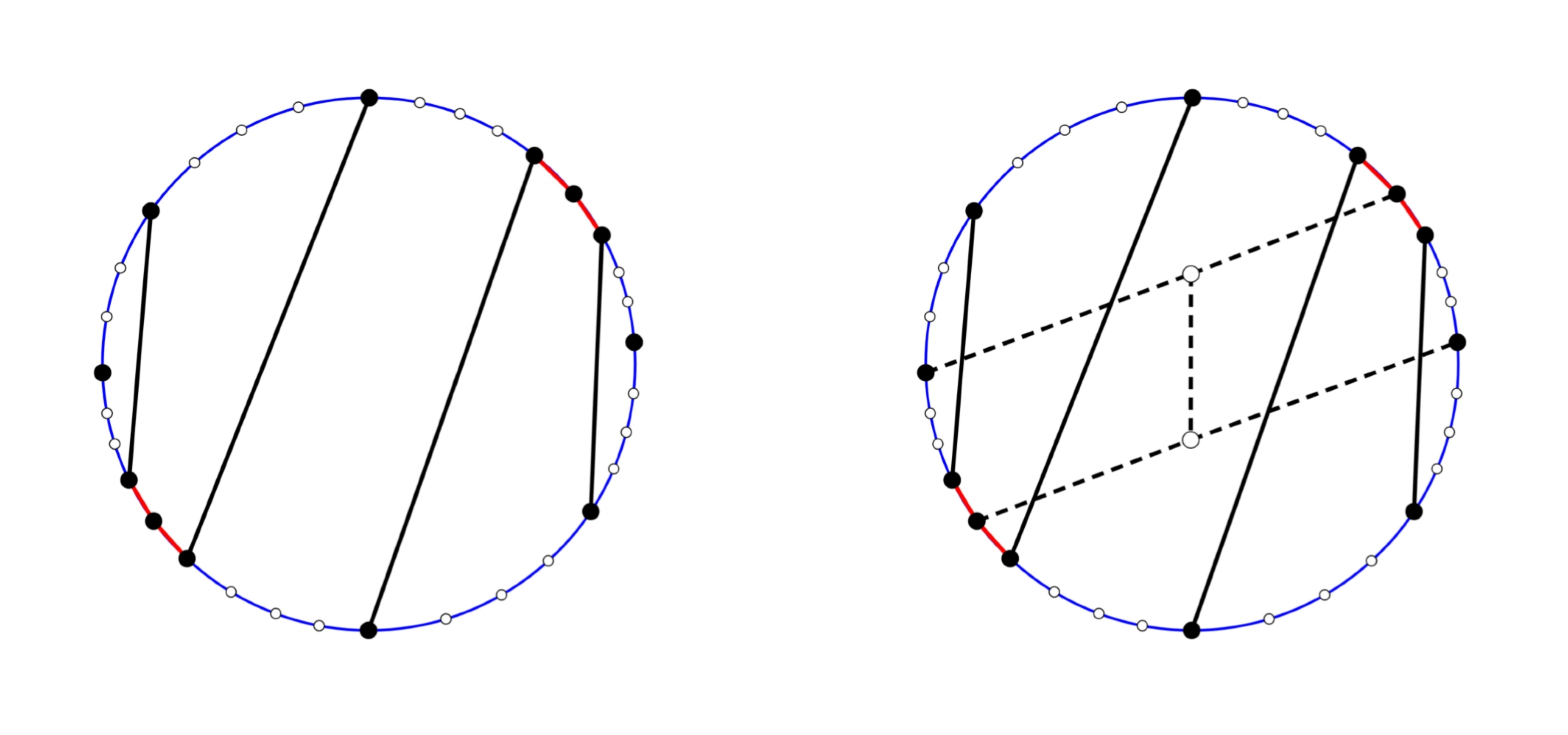}
\put(-279,155){\makebox(3,3){$z_1$}}
\put(-338,125){\makebox(3,3){$z_2$}}
\put(-350,82){\makebox(3,3){$z_3$}}
\put(-343,55){\makebox(3,3){$z_4$}}
\put(-336,40){\makebox(3,3){$z_5$}}
\put(-325,31){\makebox(3,3){$z_6$}}
\put(-279,12){\makebox(3,3){$z_7$}}
\put(-218,45){\makebox(3,3){$z_8$}}
\put(-208,89){\makebox(3,3){$z_9$}}
\put(-214,117){\makebox(3,3){$z_{10}$}}
\put(-222,129){\makebox(3,3){$z_{11}$}}
\put(-233,140){\makebox(3,3){$z_{12}$}}
\put(-88,155){\makebox(3,3){$z_1$}}
\put(-147,125){\makebox(3,3){$z_2$}}
\put(-159,82){\makebox(3,3){$z_3$}}
\put(-152,55){\makebox(3,3){$z_4$}}
\put(-145,40){\makebox(3,3){$z_5$}}
\put(-135,31){\makebox(3,3){$z_6$}}
\put(-88,12){\makebox(3,3){$z_7$}}
\put(-27,45){\makebox(3,3){$z_8$}}
\put(-17,89){\makebox(3,3){$z_9$}}
\put(-23,117){\makebox(3,3){$z_{10}$}}
\put(-31,129){\makebox(3,3){$z_{11}$}}
\put(-42,140){\makebox(3,3){$z_{12}$}}
\put(-86,116){\makebox(3,3){$m\!+\!1$}}
\put(-90,57){\makebox(3,3){$m\!+\!2$}}
\put(-279,0){\makebox(3,3){$G^{**}$}}
\put(-88,0){\makebox(3,3){$G^*$}}

\centering 
\vspace{0.1cm}
\caption {The graphs $G^{**}$ and $G^*$.}
\end{figure}
\FloatBarrier 
The following lemma shows that some pairs of vertices are $(G^*,H^*)$-exchangeable.

\begin{lem} (Alon, Defant, and Kravitz \cite{ADK}) \label{alon}
Let $G^*$ and $H^*$ be the graphs as described above. Then the two vertices $m+1$ and $m+2$ are $(G^*,H^*)$-exchangeable from the identity bijection Id $: [m+2] \mapsto [m+2]$.
\end{lem}

 The last lemma consider when the pair $(G^{**},H^{**})$ can be embedded into a pair of large random graphs. It is worth noting that  the lower bound restriction ``$p_0/\l$'' for $p_1,p_2$ is not crucial here, but it is sufficient for our use to prove Proposition \ref{connected}.
\begin{lem} \label{3}
Let $n$ be a large enough integer  , $m=\lfloor{(\log n)^{2/3}}\rfloor$, $\l=\lfloor m^{1/2}/2 \rfloor$, and  $G^{**}, H^{**}$ be two graphs as described above. Let  $\Gamma=\{x_1,\dots,x_\l,y_1,\dots,y_\l,z_3,z_5,z_9,z_{11}\}$, $q_i=\lfloor{p_0n/(5\l)}\rfloor$ for all $i\in \Gamma$, and $q_i=\lfloor n/(2m)\rfloor$ for all $i\in [m] \backslash \Gamma$.  Let $X$ and $Y$ be independently chosen random graphs in $\mathcal{G}(n,p_1)$ and $\mathcal{G}(n,p_2)$, respectively, where $p_1=p_1(n)$ and $p_2=p_2(n)$ depend on $n$, and $p_0 =\exp(2(\log n)^{2/3})/n^{1/2}$. If
\[
\begin{split}
p_1 p_2 \ge p_0^2  ~\text{ and }~  p_2\ge p_1\ge \frac{1}{{\l}}p_0,  
\end{split}
\]
 then the probability that the pair $(G^{**},H^{**})$ is $(q_1,\dots,q_m)$-embeddable in $(X,Y)$ is at least $1-n^{-n/3}$. 
\end{lem}
{\bf \pf}Assume without loss of generality that $p_1 p_2 = p_0^2$, which implies $p_2 \le \l p_0$.  We omit the floor symbols in $m,\l,q_1,\dots,q_m$ since this does not affect the asymptotic properties.  For each set $J\subseteq [m]$, let  $\gamma(J)=|J\cap\Gamma|$. Note that $Q=q_1+\dots+q_m$ satisfying $n/3 \le Q \le n$. The result will follow from Lemma \ref{a} if we can show that 
\begin{align*} 
&\, p_1^{|E(G^{**}|_J)} p_2^{|E(H^{**}|_J)|} \prod_{j\in J} q_j \\
= &\, p_1^{|E(G^{**}|_J)|} p_2^{|E(H^{**}|_J)|} {\left(\frac{p_0n}{5\l} \right)}^{\gamma(J)} {\left(\frac{n}{2m} \right)}^{|J|-\gamma(J)}  \ge 3\cdot 2^{m+1} n \log{n} 
\tag{1}
\end{align*}
 for every $J\subseteq [m]$ satisfying $|E(G^{**}|_J)|+|E(H^{**}| _J)|\ge 1$.

  If $|E(G^{**}|_J)|+|E(H^{**}| _J)|\ge 1$ and $w\notin J$, where $w$ is the center of the star  $H^{**}$, then the graph $H^{**}|_J$ has no edges, so there is an edge $\{t_1,t_2\}$ in $G^{**}|_J$.  And one of the vertices in $\{t_1,t_2\}$, say $t_1$, is not in $\Gamma$, since $\Gamma$ is an independent set in $G^{**}$. Let $J'=(J\backslash \{t_1\})\cup\{w\}$, and observe that $|J'|=|J|$, $\gamma(J')=\gamma(J)$ since $w$ is also not in $\Gamma$. There are $|J|-1$ edges in $H^{**}|_{J'}$ and at most $|J|-1$ edges in $G^{**}|_{J'}$ that are incident to $t_1$. Consequently,
\begin{align*}
&p_1^{|E(G^{**}|_J)|} p_2^{|E(H^{**}|_J)|} {\left(\frac{p_0n}{5\l} \right)}^{\gamma(J)} {\left(\frac{n}{2m} \right)}^{|J|-\gamma(J)} \\
\ge &\left(  \frac{1}{\l}\right)^{|J'|} p_1^{|E(G^{**}|_{J'})|} p_0^ {|E(H^{**}|_{J'})|} {\left(\frac{p_0n}{5\l} \right)}^{\gamma(J')} {\left(\frac{n}{2m} \right)}^{|J'|-\gamma(J')} \\
 \ge &\left(  \frac{1}{\l^2}\right)^{|J'|} p_1^{|E(G^{**}|_{J'})|} p_2^{|E(H^{**}|_{J'})|} {\left(\frac{p_0n}{5\l} \right)}^{\gamma(J')} {\left(\frac{n}{2m} \right)}^{|J'|-\gamma(J')} . \tag{2}
\end{align*}
This is to say, to show (1) holds, it suffices to prove that the right side of (2) is greater or equal to the right side of (1)  for all sets $J\subseteq [m]$ satisfying $|E(G^{**}|_J)|+|E(H^{**}| _J)|\ge 1$ and $w\in J$. Assume that $J$ satisfies these conditions, and observe that $|E(H^{**}|_J)|=|J|-1$ since $w\in J$. Let $\alpha(J)=|E(G^{**}|_J)|$. With this notation,   
\begin{align*}
&\left(  \frac{1}{\l^2}\right)^{|J|} p_1^{|E(G^{**}|_J)|} p_2^{|E(H^{**}|_J)|} {\left(\frac{p_0n}{5\l} \right)}^{\gamma(J)} {\left(\frac{n}{2m} \right)}^{|J|-\gamma(J)} \\
\ge& \left( \frac{1}{\l^5} \right)^{|J|} p_0^{\alpha(J)+|J|-1}  {\left(\frac{p_0n}{5\l} \right)}^{\gamma(J)} {\left(\frac{n}{2m} \right)}^{|J|-\gamma(J)} \\
=&\, p_0^{\alpha(J)-|J|+\gamma(J)-1}  \left(\frac{2m}{5\l}\right)^{\gamma(J)} \left( \frac{p_0^2 n}{2m\l^5}\right)^{|J|}\\ 
\ge &\, p_0^{\alpha(J)-|J|+\gamma(J)-1}   \left( \frac{p_0^2 n}{2m\l^5}\right)^{|J|}.
\end{align*}
Therefore, it suffices to prove that 
\begin{align*}
p_0^{\alpha(J)-|J|+\gamma(J)-1}   \left( \frac{p_0^2 n}{2m\l^5}\right)^{|J|}\ge 3\cdot 2^{m+1} n \log{n}. \tag{3}
\end{align*}
for all sets   $J\subseteq [m]$ satisfying $|E(G^{**}|_J)|+|E(H^{**}| _J)|\ge 1$ and $w\in J$.
 
 In what follows, let $c(J)$ be the number of components of $G^{**}|_J$. Let us also recall that $m=(\log n)^{2/3}$ and $\l=m^{1/2}/2=(\log n)^{1/3}/2$. Furthermore, $p_0^2 n/(2m\l^5)$ is certainly greater than $1$ when $n$ is large enough.

\vskip 2mm
We distinguish the following two cases according to $|J|$ separately.
\vskip 2mm
\noindent {\textbf{Case 1.}} $|J|\ge m/6$. 
\vskip 2mm
Because $G^{**}$ consists of a cycle with $4$ chords, we have $\alpha(J)\le |J|+4$, so
\[
\begin{split}
&\, p_0^{\alpha(J)-|J|+\gamma(J)-1}   \left( \frac{p_0^2 n}{2m\l^5}\right)^{|J|}\\
 \ge &\, p_0^{\gamma(J)+3} \left( \frac{p_0^2 n}{2m\l^5}\right)^{|J|}
 \ge\,  p_0^{|\Gamma|+3} \left( \frac{p_0^2 n}{2m\l^5}\right)^{m/6}\\
=& \, p_0^{2\l+7} \left( \frac{p_0^2 n}{2m\l^5}\right)^{m/6} \ge  n^{-(\l+7/2)} \left( \frac{\exp (4(\log n)^{2/3})}{2m\l^5}\right)^{m/6}\\
= & \exp \left(  -\left( \frac{1}{2}(\log n)^{1/3}+\frac{7}{2} \right) \log n \right) \left( \frac{\exp (4(\log n)^{2/3})}{(\log n)^{7/3}/16}\right)^{(\log n)^{2/3}/6}\\
=& \exp \left( \frac{1}{6}(\log n)^{4/3}+ O(\log n) \right) \ge 3\cdot 2^{m+1} n \log{n}. \\
\end{split}
\] 
\vskip 2mm
\noindent {\textbf{Case 2.} $|J|<m/6$.  
\vskip 2mm
 The graph $G^{**}$ has girth at least $m/6$, which forces the induced subgraph $G^{**}|_J$ to be a forest. So we have $\alpha(J)=|J|-c(J)$ and we can rewrite (3) as 
\begin{align*}
p_0^{\gamma(J)-c(J)-1}   \left( \frac{p_0^2 n}{2m\l^5}\right)^{|J|}\ge 3\cdot 2^{m+1} n \log{n}. 
\end{align*}
\vskip 2mm
\noindent {\textbf{Subcase 2.1.}} $\gamma(J)\le c(J)-1$. 
\vskip 2mm

Note that $|J|\ge 2$ since $|E(G^{**}|_J)|+|E(H^{**}| _J)|\ge 1$. Consequently,
\[
\begin{split}
  p_0^{\gamma(J)-c(J)-1} \left( \frac{p_0^2 n}{2m\l^5}\right)^{|J|} & \ge  p_0^{-2}\left( \frac{p_0^2 n}{2m\l^5}\right)^{2} = \frac{ (p_0 n)^2}{(2m\l^5)^2} \\&= \frac{\exp \left(4(\log n)^{2/3} \right)n}{(2m\l^5)^2}  \ge 3\cdot 2^{m+1} n \log{n}. \\
\end{split}
\] 
\vskip 2mm
\noindent {\textbf{Subcase 2.2.}} $c(J)\le \gamma(J) \le c(J)+1$. 
\vskip 2mm

The number of elements of $\Gamma \cup \{w\}$ that are in $J$ is $\gamma(J)+1\ge c(J)+1$. This means that some component of $G^{**}|_J$ contains at least $2$ elements of $\Gamma \cup \{w\}$. The minimum distance in $G^{**}$ between any two elements of $\Gamma \cup \{w\}$ is $\l-1$, so $|J|\ge \l$. It follows that
\[
\begin{split}
p_0^{\gamma(J)-c(J)-1}   \left( \frac{p_0^2 n}{2m\l^5}\right)^{|J|} &\ge  \left( \frac{p_0^2 n}{2m\l^5}\right)^{\l} \ge  \frac{\exp \left( 4\l(\log n)^{2/3} \right)}{(2m\l^5)^{\l}}  \\
&=n^{2+o(1)}
 \ge  3\cdot 2^{m+1} n \log{n}. \\
\end{split}
\] 
\vskip 2mm
\noindent {\textbf{Subcase 2.3.}} $ \gamma(J)\ge c(J)+2$.
\vskip 2mm

Recall that for any distinct $s_1,s_2\in \Gamma \cup \{w\}$, if $\{s_1,s_2\}\neq \{z_3,z_5\}$ and $\{s_1,s_2\}\neq \{z_9,z_{11}\}$, then the distance between $s_1$ and $s_2$ in $G^{**}$ is at least $m/(3\l)$. Since $w\in J$, we have $|J|\ge (\gamma(J)-c(J)-1)m/(3\l)$. Indeed, if a component of $G^{**}|_J$ contains $k$ elements of $(\Gamma \cup \{w\})\backslash \{z_3,z_9\}$, then this component must contain at least $(k-1)m/(3\l)$ vertices, and $J$ contains precisely $\gamma(J)-1$ elements of  $(\Gamma \cup \{w\})\backslash \{z_3,z_9\}$. Therefore
\[
\begin{split}
p_0^{\gamma(J)-c(J)-1}   \left( \frac{p_0^2 n}{2m\l^5}\right)^{|J|} &\ge \left( n^{-1/2}\frac{(p_0^2n)^{m/(3\l)}}{(2m\l^5)^{m/(3\l)}} \right)^{\gamma(J)-c(J)-1}\\
&= \left( n^{-1/2+o(1)}(p_0^2 n)^{m/(3\l)} \right)^{\gamma(J)-c(J)-1} \\
&\ge \left(
 n^{-1/2+o(1)}\left(
 \exp \left( 4(\log n)^{2/3} \right)  
\right)^{(2/3)(\log n)^{1/3}} 
\right)^{\gamma(J)-c(J)-1} \\
&=  \left( n^{-1/2+8/3+o(1)} \right)^{\gamma(J)-c(J)-1} \ge n^{13/6+o(1)}\\
&\ge 3\cdot 2^{m+1} n \log{n}. \\
\end{split}
\] 
The proof of Lemma \ref{3} is complete.\qed

\vskip 3mm
We now begin to prove Proposition \ref{connected}.
\vskip 2mm
\noindent {\bf \emph{Proof of Proposition} \ref{connected}}.   By Lemma \ref{dk}, we may assume $p_1\le p_2$. Let $n$ be large enough, and $m$, $G^*$, $G^{**}$, $ H^*$ and $H^{**}$ be as described above. Let $\Gamma$ and $q_1,\dots,q_m$ be as in the statement of Lemma \ref{3}. Thus,  the pair $(G^{**},H^{**})$ is $(q_1,\dots,q_m)$-embeddable in $(X,Y)$ w.h.p. by Lemma \ref{3}. Because both $p_1$ and $p_2$ are  much larger than $\log{n} /n$, it is well known that  both $X$ and $Y$ are connected w.h.p. and the degrees of all vertices in $X$ and $Y$ are $p_1 n(1+o(1))$ and $p_2 n(1+o(1))$ w.h.p., respectively. Hence, we may assume that $X$ and $Y$ have these properties.

 Choose arbitrarily vertices $u,v \in V(Y)$ and a bijection $\sigma : V(X)\mapsto V(Y)$ such that $\{ \sigma^{-1}(u),\sigma^{-1}(v)\} \in E(X)$. By Lemma \ref{ex}, it suffices to show that $u,v$ are $(X,Y)$-exchangeable from $\sigma$. Let us choose pairwise disjoint subsets $V_1,\dots,V_m$ of $V(Y)\backslash \{u,v\}$ such that
\begin{itemize}
\item  $|V_i|=q_i$ for all $i\in [m]$;
\item $\sigma^{-1}(V_{z_3})$ and $\sigma^{-1}(V_{z_{11}})$ are contained in the neighborhood of $\sigma^{-1}(u)$ in $ X$;
\item $\sigma^{-1}(V_{z_5})$ and $\sigma^{-1}(V_{z_{9}})$ are contained in the neighborhood of $\sigma^{-1}(v)$ in $ X$;
\item  $V_{x_1},\dots,V_{x_\ell}$ are all contained in the neighborhood of $u$ in $Y$;
\item  $V_{y_1},\dots,V_{y_\ell}$ are all contained in the neighborhood of $v$ in $Y$.
\end{itemize}

  Such a choice is possible since the condition $q_i=\lfloor p_0 n/ (5\l) \rfloor$
 for all $i\in \Gamma$ guarantees that $4q_i \le p_1 n (1+o(1))$ and $(2\l+4)q_i \le p_2 n (1+o(1))$  for all $i\in \Gamma$. Because the pair $(G^{**},H^{**})$ is $(q_1,\dots,q_m)$-embeddable in $(X,Y)$, it must be the case that $(G^{**},H^{**})$ is embeddable in $(X,Y)$ with respect to the sets $V_1,\dots,V_m$ and the bijection $\sigma$. This means that there exist vertices $v_i\in V_i$ for all $i\in[m]$, such that for all $i,j\in[m]$, we have
\[
\begin{split}
\{i,j\} \in E(H) \Rightarrow \{v_i,v_j\} \in E(Y)  \text{ and } \\
\{i,j\} \in E(G) \Rightarrow \{ \sigma^{-1}(v_i),\sigma^{-1}(v_j)\} \in E(X).
\end{split}
\]  

  Define a map $\psi: V(H^*) \mapsto V(Y)$ by $\psi(m+1)=u$, and $\psi(m+2)=v$, and $\psi(i)=v_i$ for all $i\in [m]$. Define $\varphi: V(G^*) \mapsto V(X)$ by $\varphi =\sigma^{-1}\circ \psi \circ \text{Id}$. It is immediate from our construction that $\psi$ is a graph embedding of $H^*$ into $Y$ satisfying $\psi(m+1)=u$ and $\psi(m+2)=v$. Similarly,  $\varphi$ is  a  graph embedding  of  $G^*$ into $X$. By
  Lemma \ref{alon}, we know  that the vertices $m+1$ and $m+2$ are $(G^*,H^*)$-exchangeable from $Id$ $: [m+2] \mapsto [m+2]$. By Lemma \ref{ex2},  $u,v$ are $(X,Y)$-exchangeable from $\sigma$.  
 
 The proof of Proposition \ref{connected} is complete. \qed 
 
 \vskip 3mm
\noindent{\bf Concluding Remark.} Combining Propositions \ref{disconnected} and \ref{connected}, we complete the proof of Theorem \ref{mr}. From the proof of Proposition \ref{connected},  one can see that the ``$1/\l$" term in the lower bound $(1/\l)p_0$ in Theorem \ref{mr} comes from roughly the ratio of the number of vertices that are adjacent to the vertex $m+1$ (or the vertex $m+2$) in $G^{**}$ and $H^{**}$. Along this way, it seems hard to improve the term ``$1/\l$'' to ``$1/n^c$'' for some constant $c>0$, unless we can find some suitable pair of small sparse graphs other than the pair $G^{**}$ and $H^{**}$.

\section*{Acknowledgments}

 We are grateful to the anonymous referees for their very careful comments. This research was supported by NSFC under grant numbers 11871270, 12161141003 and 11931006.

\end{document}